\RequirePackage{fix-cm}
\documentclass{svjour3}                     
\smartqed  
\usepackage{graphicx}
\usepackage{amsmath}
\usepackage{amssymb}
\usepackage{epsfig}
\usepackage{mathrsfs}
\usepackage{amsfonts}
\def\pn{\par\smallskip\noindent}
\def\proof{\pn {Proof.} }
\def\endproof{\hfill \quad{$\Box$}\smallskip}

\begin{document}

\title{On Local Convexity of Quadratic Transformations
\thanks{This research was supported by Beijing Higher Education Young Elite Teacher Project 29201442, and by the fund of State Key Laboratory of Software Development Environment under grant SKLSDE-2013ZX-13.}
}

\titlerunning{On Local Convexity of Quadratic Transformations}        

\author{Yong Xia}


\institute{Y. Xia \at
              State Key Laboratory of Software Development
              Environment, LMIB of the Ministry of Education, School of
Mathematics and System Sciences, Beihang University,  Beijing
100191,  P. R. China
              \email{dearyxia@gmail.com }
}

\date{Received: date / Accepted: date}

\maketitle

\begin{abstract}
In this paper, we improve Polyak's local convexity result for quadratic transformations. Extension and open problems are also presented.
\end{abstract}

\section{Introduction}
\label{intro}

Let $x\in\Bbb R^n$ and $f(x) = (f_1(x),\ldots, f_m(x))$,
where
\[
f_i(x) =
\frac{1}{2}x^TA_ix + a_i^Tx, ~i = 1,\ldots,m
\]
are quadratic functions. One interesting question is when the following joint numerical range
\[
F_m=\{f(x):x\in\Bbb R^n\}\subseteq \Bbb R^m
\]
is convex.

The first such result is due to Dines \cite{Di} in 1941.
It states that if $f_1, f_2$ are homogeneous quadratic
functions then the set $F_2$
is convex. In 1971, Yakubovich \cite{y71,y77} used this basic result to prove the famous S-lemma, see \cite{P} for a survey.
Brickman \cite{Br} proved in 1961 that if $f_1, f_2$ are homogeneous quadratic
functions and $n\ge 3$ then the set
 $\{(f_1(x), f_2(x)) :  x\in \Bbb R^n, \|x\| = 1\} \subseteq \Bbb R^2$
is convex.
Fradkov \cite{Fr} proved in 1973 that if matrices $A_1, \ldots, A_m$ commute and $f_1,\ldots,f_m$ are homogeneous, then $F_m$ is convex.
In 1995, it was showed by Ramana and Goldman \cite{Ra} that the identification of
the convexity of $F_m$ is NP-hard. In the same paper, the quadratic maps, under which the image of
every linear subspace is convex, was also investigated.
Based on Brickman's result, Polyak \cite{P98} proved in 1998 that
if $n\ge 3$ and $f_1, f_2, f_3$ are homogeneous
quadratic functions such that $\mu_1A_1+\mu_2A_2+\mu_2A_3\succ 0$ (where notation $A\succ 0$ means that $A$ is positive definite) for some $\mu\in \Bbb R^3$, then the set
 $F_3$ is convex. Moreover, as shown in the same paper, when $n\ge 2$ and there exists $\mu\in\Bbb R^2$ such that $\mu_1A_1+\mu_2A_2\succ 0$, the set
 $F_2$ is convex.
In 2007, Beck \cite{Be} showed that  if $m\le n$,  $A_1\succ 0$ and $A_2=\ldots=A_m=0$, then $F_m$ is convex. However,  if $A_1\succ 0$, $A_2=\ldots=A_{n+1}=0$ and $a_2,\ldots,a_{n+1}$ are linearly independent, then $F_{n+1}$ is not convex. When $m=2$, Beck's result reduces to be a corollary of Polyak's result. Very recently, Xia et al. \cite{Xia} used the new developed S-lemma with equality to establish the necessary and sufficient condition for the convexity of $F_2$ for $A_2=0$ and arbitrary $A_1$.

More generally, Polyak \cite{P01,P03} succeeded in proving a nonlinear image of a small ball in a Hilbert space is convex, provided that the map is
$C^{1,1}$
and the center of the ball is a regular point of the map. Later, Uderzo \cite{U} extended the result to a certain subclass of uniformly convex Banach spaces.
When focusing on quadratic transformations, Polyak's result reads as follows:
\begin{theorem}[\cite{P}]\label{thm:1}
Let $A=[a_1~\ldots~a_m]\in \Bbb R^{n\times m}$  and define
\begin{eqnarray}
&&L:=\sqrt{\sum_{i=1}^m\|A_i\|^2},\label{L} \\
&&\nu:=\sigma_{\min}(A)=\sqrt{\lambda_{\min}(A^TA)},\nonumber
\end{eqnarray}
where $\|A_i\|=\sigma_{\max}(A_i)=\sqrt{\lambda_{\max}(A_i^TA_i)}$ is the spectral
norm of $A_i$, $\sigma_{\min}(\cdot)$, $\lambda_{\min}(\cdot)$, $\sigma_{\max}(\cdot)$, $\lambda_{\max}(\cdot)$, denote the smallest and largest singular
value and eigenvalue, respectively.

If $\epsilon< \epsilon^*:= \nu/(2L)$,
then the image
\begin{equation}
F_m(\epsilon) = \{f(x):~x\in\Bbb R^n,~\|x\|\le \epsilon\}\label{F}
\end{equation}
is a convex set in $\Bbb R^m$.
\end{theorem}
Polyak \cite{P01,P03} used  the following example to show his estimation $\epsilon^*$ is tight,
where $n=m=2$ and
\[
f_1(x)=x_1x_2-x_1,~f_2(x)=x_1x_2+x_2.
\]
Actually, in this case, $\epsilon^*=1/(2\sqrt{2})
\approx0.3536$. It is trivially verified that $F_m(\epsilon)$ is convex for
$\epsilon\le \epsilon^*$ and loses
convexity for $\epsilon>\epsilon^*$.

In this paper, we improve the above Polyak's result for quadratic transformations (i.e., Theorem \ref{thm:1}) by strengthening the constant $L$. Then, Theorem \ref{thm:1} is extended to the image of the ball of the same  radius $\epsilon$ centered at any point $a$ satisfying $\|a\|<2(\epsilon^*-\epsilon)$. Furthermore, we propose two new approaches
for possible improvement of $L$.

The paper is organized as follows. In Section 1, we improve and extend Theorem \ref{thm:1}.
In Section 2, we discuss further possible improvements. In the final conclusion section, we propose two open questions.

Throughout the paper, all vectors are column vectors.
Let $v(\cdot)$ denote the optimal value of
problem $(\cdot)$.  Notation $A\succeq 0$ implies that
the matrix $A$ is positive semidefinite. ${\rm vec}(A)$ denotes the vector obtained by stacking the columns of $A$ one underneath the
other. The trace of $A$ is denoted by trace$(A)=\sum_{i=1}^nA_{ii}$.
The
Kronecker product and the inner product of the matrices $A$ and $B$ are denoted by
$A\otimes B$ and $A\bullet B={\rm
trace}(AB^T)=\sum_{i,j=1}^na_{ij}b_{ij}$, respectively. The
identity matrix is denoted by $I$. $\|x\|=\sqrt{x^Tx}$ is the standard norm of the vector $x$.

\section{Main Results}
In this section, we first improve Theorem \ref{thm:1} and then extend it to the ball of the same  radius  centered at any point close enough to the zero point.
\begin{theorem}\label{thm:2}
Define
\begin{equation}
 L_{\rm new}:=\sqrt{\lambda_{\max}\left(\sum_{i=1}^mA_i^TA_i\right)}.
 \label{Lnew}
\end{equation}
Then we have
\begin{equation}
L_{\rm new} \le L.\label{LL}
\end{equation}
For any $\epsilon< \epsilon^*_{\rm new}:= \nu/(2L_{\rm new})$,
the image $F_m(\epsilon)$ defined in (\ref{F})
is convex.
\end{theorem}
\proof
Let $L_b$ be
any upper bound of the Lipschitz constant of $f$, i.e.,
\begin{equation}
\|\nabla f(x)-\nabla f(z)\|\le L_b\|x-z\|,~ \forall x, z\in \Bbb R^n.\label{Lb}
\end{equation}
According to the proof in \cite{P01},
Theorem \ref{thm:1} remains true if
$L$ defined in (\ref{L}) is replaced by $L_b$.
It is sufficient to show that $L_b:=L_{\rm new}$ satisfies (\ref{Lb}). To this end, we have
\begin{eqnarray}
&&\max_{\|x-z\|=1}\|\nabla f(x)-\nabla f(z)\|\nonumber\\
&=&\max_{\|x-z\|=1}\|[A_1(x-z)~\ldots~A_m(x-z)]\|\nonumber\\
&=&\max_{\|y\|=1}\|[A_1y~\ldots~A_my]\|\nonumber\\
&=&\sqrt{\max_{\|y\|=1} \lambda_{\max}\left([A_1y~\ldots~A_my]^T[A_1y~\ldots~A_my]\right)} \label{ineq0}\\
&\le&\sqrt{\max_{\|y\|=1}  {\rm trace }\left([A_1y~\ldots~A_my]^T[A_1y~\ldots~A_my]\right)}\label{ineq1}\\
&=&\sqrt{\max_{\|y\|=1} y^T\left(\sum_{i=1}^mA_i^TA_i\right)y}\nonumber\\
&=&\sqrt{\lambda_{\max}\left(\sum_{i=1}^mA_i^TA_i\right) }.\nonumber
\end{eqnarray}
The inequality (\ref{LL}) holds since
\begin{eqnarray*}
L_{\rm new} &=& \sqrt{\lambda_{\max}\left(\sum_{i=1}^mA_i^TA_i\right)}
=
\sqrt{\max_{\|y\|=1}y^T\left(\sum_{i=1}^mA_i^TA_i\right)y}\\
&\le&
\sqrt{ \sum_{i=1}^m \left(\max_{\|y\|=1}y^TA_i^TA_iy\right) }
=
\sqrt{\sum_{i=1}^m\lambda_{\max}\left(A_i^TA_i\right)}
=
\sqrt{\sum_{i=1}^m\|A_i\|^2}=L.
\end{eqnarray*}
\endproof

\begin{theorem}\label{thm:3}
For any $0<\epsilon< \epsilon^*_{\rm new}=\nu/(2L_{\rm new})$ and any $a\in\Bbb R^n$ such that $\|a\|<2(\epsilon^*_{\rm new}-\epsilon)$,
the image
\[
F_m( \epsilon,a)=\{f(x):~x\in\Bbb R^n,~\|x-a\|\le \epsilon \}
\]
is a convex set in $\Bbb R^m$.
\end{theorem}
\proof
For any $a\in\Bbb R^n$ such that $\|a\|<2(\epsilon^*_{\rm new}-\epsilon)$, we have
\begin{eqnarray}
&&\sigma_{\min}(A+[A_1a\ldots A_ma]) \nonumber \\
&\ge& \sigma_{\min}(A)-\sigma_{\max}(-[A_1a\ldots A_ma])\nonumber\\
&\ge& \sigma_{\min}(A)-\sup_{\|a\|<2(\epsilon^*_{\rm new}-\epsilon)}\sigma_{\max}([A_1a\ldots A_ma])\nonumber\\
&=&\sigma_{\min}(A)- \sqrt{\sup_{\|a\|<2(\epsilon^*_{\rm new}-\epsilon)} \lambda_{\max}\left([A_1a~\ldots~A_ma]^T[A_1a~\ldots~A_ma]\right)} \nonumber\\
&\ge&\sigma_{\min}(A)- \sqrt{\sup_{\|a\|<2(\epsilon^*_{\rm new}-\epsilon)}   {\rm trace }\left([A_1a~\ldots~A_ma]^T[A_1a~\ldots~A_ma]\right)} \nonumber\\
&=&\sigma_{\min}(A)- \sqrt{\sup_{\|a\|<2(\epsilon^*_{\rm new}-\epsilon)}
a^T\left(\sum_{i=1}^mA_i^TA_i\right)a}\label{ieq}\\
&=&\sigma_{\min}(A)-2(\epsilon^*_{\rm new}-\epsilon) \sqrt{\lambda_{\max}\left(\sum_{i=1}^mA_i^TA_i\right) }\nonumber\\
&=&\sigma_{\min}(A)-2(\epsilon^*_{\rm new}-\epsilon)L_{\rm new}\nonumber\\
&=& 2\epsilon L_{\rm new},\nonumber
\end{eqnarray}
where the first inequality is Weyl's inequality \cite{H} for the singular values, see also
 Problem III.6.5 in \cite{B} or  Theorem 3.3.16 in \cite{H2}.

Since the optimal value of the maximizing problem (\ref{ieq}) is unattainable, the above inequality implies that
\begin{equation}
\sigma_{\min}(A+[A_1a\ldots A_ma])>2\epsilon L_{\rm new},~\forall a\in\Bbb R^n:~\|a\|<2(\epsilon^*_{\rm new}-\epsilon).\label{ieq2}
\end{equation}
Notice that
\[
f_i(x) =f_i(a)+(A_ia+a_i)^T(x-a)+
\frac{1}{2}(x-a)^TA_i(x-a), ~i = 1,\ldots,m.
\]
Then, we have
\[
F_m( \epsilon,a)-f(a)
= \{g(y):~y\in\Bbb R^n,~\|y\|\le \epsilon \}:=G_m(\epsilon,a),
\]
where  $g(y) = ((A_1a+a_1)^Ty+
\frac{1}{2}y^TA_1y,\ldots, (A_ma+a_m)^Ty+
\frac{1}{2}y^TA_my)$.
According to Theorem \ref{thm:2},
for any
\begin{equation}
\epsilon<  \sigma_{\min}(A+[A_1a\ldots A_ma])/(2L_{\rm new}), \label{ieq3}
\end{equation}
the image $G_m(\epsilon,a)$
is a convex set in $\Bbb R^m$.
The proof is complete as (\ref{ieq3}) is ensured by (\ref{ieq2}).
\endproof

\begin{remark}
Theorem \ref{thm:2} is a special case of Theorem \ref{thm:3} by setting $a=0$.
\end{remark}

\section{Discussion}
The estimation of Theorem \ref{thm:2} is still not tight. Actually, $L_{\rm new}$ defined in (\ref{Lnew}) can be further improved to be the Lipschitz constant of $f$, denoted by $L_f$. According to
(\ref{ineq0}), we have
\begin{equation}
L_f^2=\max_{\|y\|=1} \lambda_{\max}\left([A_1y~\ldots~A_my]^T[A_1y~\ldots~A_my]\right).
\label{ineq2}
\end{equation}
However, this is a nonlinear eigenvalue optimization problem and not easy to solve.
Except for the upper bound  $L_{\rm new}$ (\ref{Lnew}), we further consider the other two relaxations of (\ref{ineq2}). We first need two lemmas.
\begin{lemma}[\cite{B}]\label{lem:1}
 Every eigenvalue of $B\in\Bbb R^{m\times m}$ lies within at least one of the Gershgorin discs
\[
\left\{\lambda:~\left|\lambda-B_{ii} \right|\le\sum_{j\neq i}|B_{ij}| \right\},~i=1,\ldots,m.
\]
\end{lemma}
\begin{lemma}[\cite{Gu}]\label{lem:2}
For any $m\times m$ matrix $B$, all its eigenvalues are located in the same disk
\begin{equation}
\left|\lambda-\frac{{\rm trace }(B)}{m}\right|\le \sqrt{\frac{m-1}{m}\left({\rm trace }(B^TB)-\frac{\left({\rm trace }(B)\right)^2}{m}\right)}.\label{lam}
\end{equation}
\end{lemma}
\begin{remark}
Let $\lambda_i(B)$ be the $i$-th largest eigenvalue of $B$.
When $B\succeq 0$,
substituting the following inequality
\[
{\rm trace }(B^TB)=\sum_{i=1}^m \lambda_i^2(B)
\le \left(\sum_{i=1}^m \lambda_i(B)\right)^2=\left({\rm trace }(B)\right)^2
\]
into (\ref{lam}),
we see that
Lemma \ref{lem:2} improves the inequality
\[
\lambda_{\max}(B)\le {\rm trace }(B),
\]
which is used in (\ref{ineq1}).
\end{remark}
Now, we apply Lemmas \ref{lem:1} and \ref{lem:2} to establish two new relaxations of $L_f$ (\ref{ineq2}).

Firstly, according to Lemma \ref{lem:1}, we have:
\begin{eqnarray*}
&&\sqrt{\max_{\|y\|=1} \lambda_{\max}\left([A_1y~\ldots~A_my]^T[A_1y~\ldots~A_my]\right)}\\
&\le&\sqrt{\max_{\|y\|=1} \max_{i=1,\ldots,m} \left\{y^T(A_i^TA_i)y+\sum_{j\neq i}y^T|A_i^TA_j|y\right\}}\\
&=&\sqrt{ \max_{i=1,\ldots,m} \max_{\|y\|=1} y^T\left(A_i^TA_i+\sum_{j\neq i}|A_i^TA_j|\right)y }\\
&=&\sqrt{ \max_{i=1,\ldots,m} \lambda_{\max}\left(A_i^TA_i+\frac{1}{2}\sum_{j\neq i}\left(|A_i^TA_j|+|A_j^TA_i|\right)\right)}\\
&:=& \overline{L}_{\rm new}.
\end{eqnarray*}
Consequently, Theorem \ref{thm:1} holds true if we replace $L$ with $\overline{L}_{\rm new}$.

Secondly,
according to Lemma \ref{lem:2},  we have:
\begin{eqnarray*}
&& \lambda_{\max}\left([A_1y~\ldots~A_my]^T[A_1y~\ldots~A_my]\right)\\
&\le&  \frac{1}{m}\sqrt{ \left(y^T\left(\sum_{i=1}^mA_i^TA_i\right)y\right)^2}+\\
&&\sqrt{\frac{m-1}{m}\left( \sum_{i,j=1}^m(y^TA_i^TA_jy)^2-
\frac{1}{m} \left(y^T\left(\sum_{i=1}^mA_i^TA_i\right)y\right)^2\right)} \\
&=&   \frac{1}{m}\sqrt{ z^T
\left(\left(\sum_{i=1}^m  A_i^TA_i\right)\otimes \left(\sum_{i=1}^m  A_i^TA_i\right)\right) z}+
\sqrt{\frac{m-1}{m}}\cdot\\
&&\sqrt{z^T\left( \sum_{i,j=1}^m (A_i^TA_j)\otimes(A_i^TA_j) -
\frac{1}{m} \left(\sum_{i=1}^m  A_i^TA_i\right)\otimes \left(\sum_{i=1}^m  A_i^TA_i\right)  \right)z}\\
&=&\frac{1}{m}\sqrt{ \left(\left(\sum_{i=1}^m  A_i^TA_i\right)\otimes \left(\sum_{i=1}^m  A_i^TA_i\right)\right)\bullet Z}+
\sqrt{\frac{m-1}{m}}\cdot\\
&&~~\sqrt{\left( \sum_{i,j=1}^m (A_i^TA_j)\otimes(A_i^TA_j) -
\frac{1}{m} \left(\sum_{i=1}^m  A_i^TA_i\right)\otimes \left(\sum_{i=1}^m  A_i^TA_i\right)  \right) \bullet Z}\\
&:=&B(Z)
\end{eqnarray*}
where $z=y\otimes y$ and $Z=zz^T$.  Since $y^Ty=1$,
we have
\begin{eqnarray*}
&&{\rm trace}(Z)=z^Tz=(y\otimes y)^T(y\otimes y)=(y^Ty)\otimes (y^Ty)=1\otimes1=1,\\
&&{\rm vec}(I)^TZ{\rm vec}(I)=\left({\rm vec}(I)^Tz\right)^2=\left(\sum_{i=1}^my_i^2\right)^2 =1,\\
&&\|Z{\rm vec}(I)\|=\|zz^T{\rm vec}(I)\|=\left|{\rm vec}(I)^Tz\right|\|z\|
=\|z\|=\sqrt{z^Tz}= 1,\\
&&Z=zz^T\succeq 0.
\end{eqnarray*}
Therefore, Theorem \ref{thm:1} remains true if $L$ is replaced by $\widetilde{L}_{\rm new}$, where
\begin{eqnarray*}
\widetilde{L}_{\rm new}^2=&\max&~B(Z)\\
&{\rm s.t.}& {\rm trace}(Z)=1,\\
&& {\rm vec}(I)^TZ{\rm vec}(I)=1,\\
&&\|Z{\rm vec}(I)\|\le 1,\\
&&Z\succeq 0,
\end{eqnarray*}
which is a convex semidefinite programming (CSDP) problem, and hence can be efficiently solved. In the following examples, the CSDP problems are modeled by CVX 1.2 \cite{Gr} and solved by
SDPT3 \cite{To} within CVX.

\begin{example}
Let $n =3,~ m = 2$.
Consider the two examples:
\begin{eqnarray*}
&(E_1):&A_1=\left[\begin{array}{ccc}2&0&6\\0&0&6\\6&6&2 \end{array}\right],~
A_2=\left[\begin{array}{ccc}6&5&2\\5&4&0\\2&0&0 \end{array}\right],~
A=\left[\begin{array}{cc}-1&0\\0&1\\0&0 \end{array}\right],\\
&(E_2):&A_1=\left[\begin{array}{ccc}0&5&3\\5&0&6\\3&6&4 \end{array}\right],~
A_2=\left[\begin{array}{ccc}0&4&2\\4&0&4\\2&4&4 \end{array}\right],~
A=\left[\begin{array}{cc}-1&0\\0&1\\0&0 \end{array}\right].
\end{eqnarray*}
We can verify that
\begin{eqnarray*}
&(E_1):&L\approx 14.4166,~ L_{\rm new}\approx13.9094,~
\overline{L}_{\rm new}\approx12.8849,~\widetilde{L}_{\rm new}\approx 12.6747,\\
&(E_2):&L\approx 13.8065,~ L_{\rm new}\approx13.8043,~
\overline{L}_{\rm new}\approx14.5901,~\widetilde{L}_{\rm new}\approx 13.8009.
\end{eqnarray*}
It is observed that
neither $L_{\rm new}$ nor $\overline{L}_{\rm new}$ dominates each other. Moreover, both are dominated by $\widetilde{L}_{\rm new}$.

Figure \ref{fig}  shows the images of the $\epsilon$-discs for $(E_1)$ and $(E_2)$, respectively.
It follows that $\widetilde{L}_{\rm new}$ is not tight and
the convexity loses when $\epsilon$ is large enough.
\begin{figure}
\begin{center}
\includegraphics[width=0.5\textwidth]{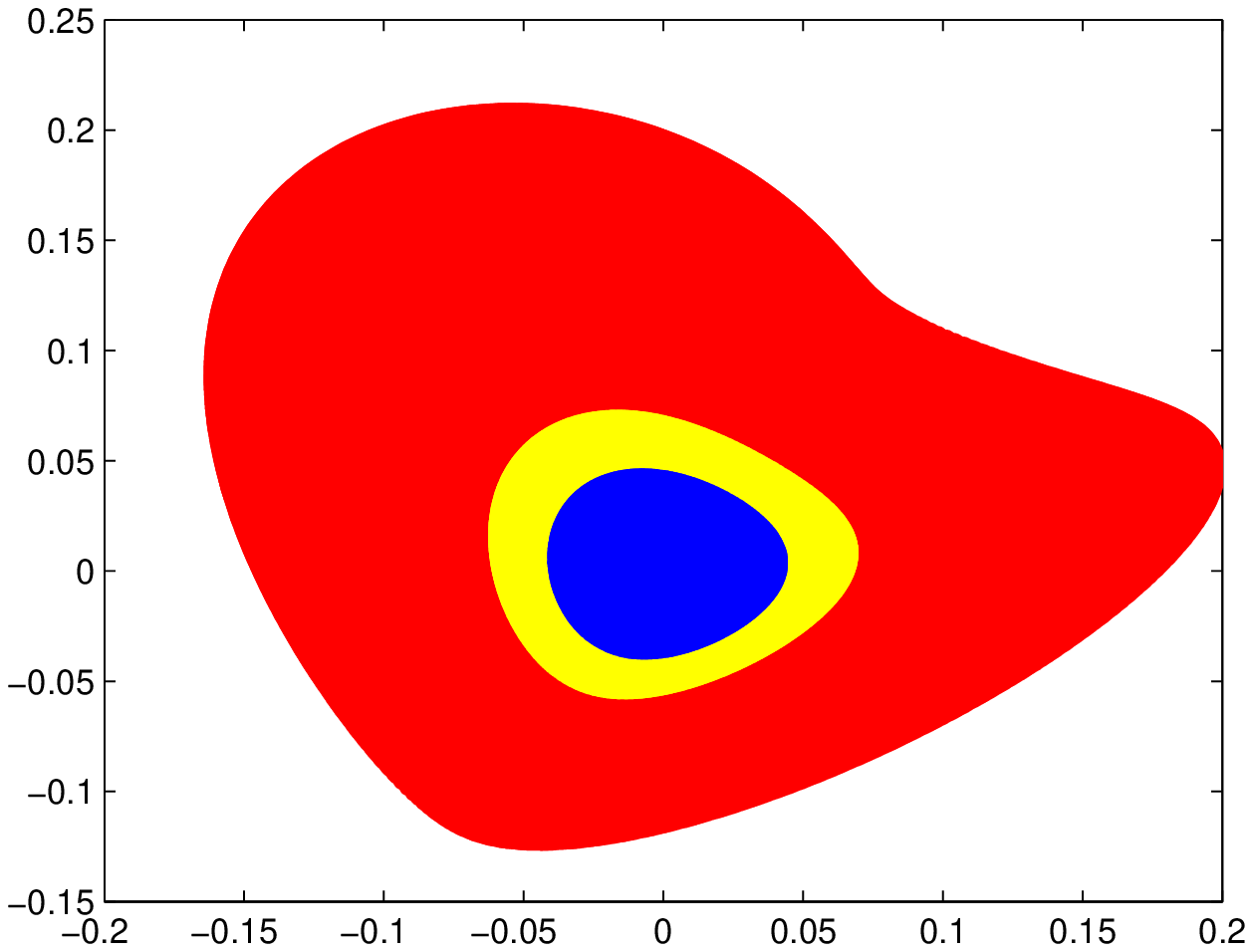}\includegraphics[width=0.5\textwidth]{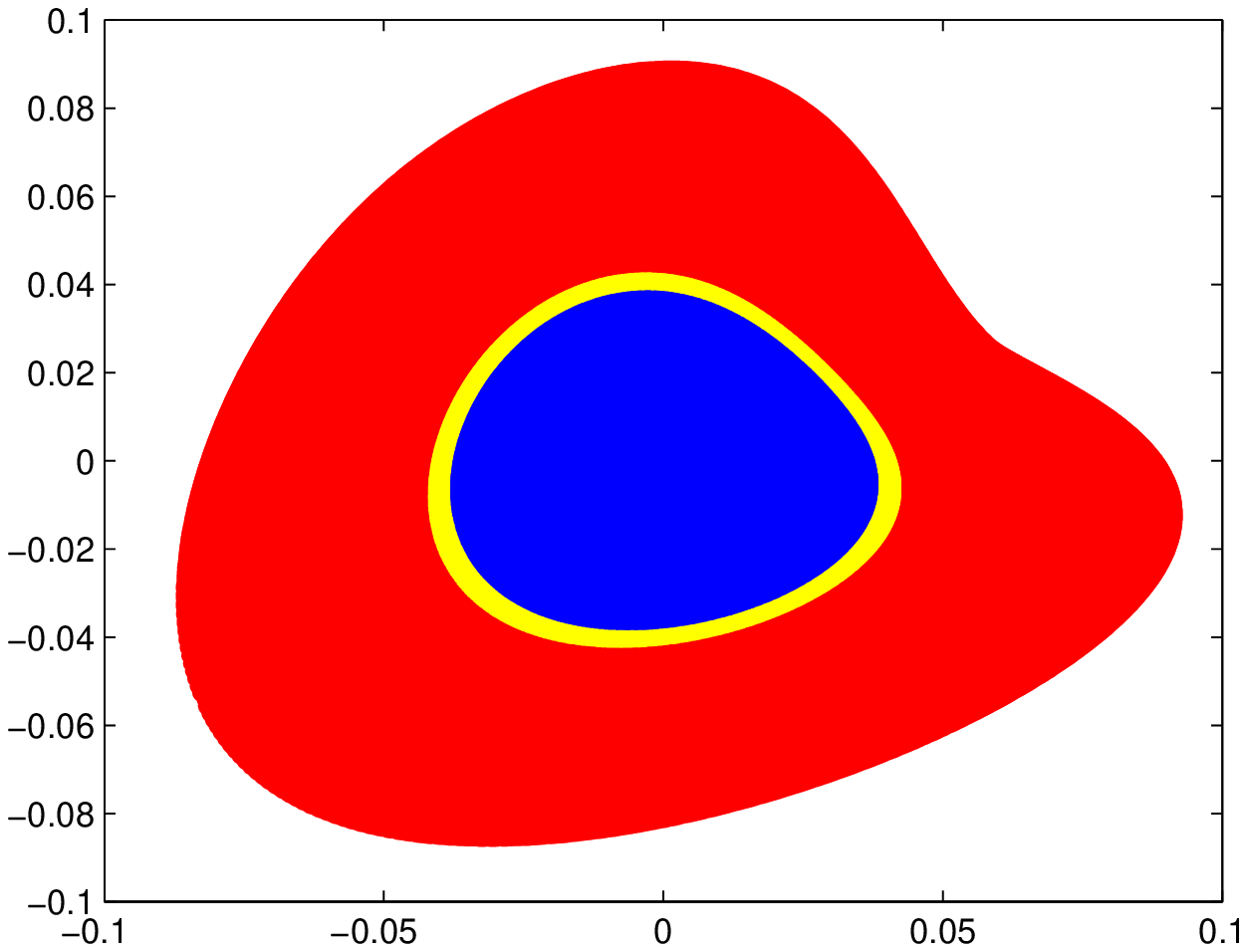}
\caption{Images of $\epsilon$-discs for $(E_1)$ with $\epsilon=1/(2\widetilde{L}_{\rm new})\approx0.0394,~0.06,~0.14$ in the left subgraph and
for $(E_2)$ with
$\epsilon=1/(2\widetilde{L}_{\rm new})\approx0.0362,~0.04,~0.08$ in the right subgraph. 
}\label{fig}
\end{center}
\end{figure}
\end{example}

\section{Conclusions}
In this paper, we improve and extend Polyak's local convexity result for quadratic transformations by providing tighter bounds for
\[
\max_{\|y\|=1} \lambda_{\max}\left([A_1y~\ldots~A_my]^T[A_1y~\ldots~A_my]\right).
\]
It is open whether the above nonlinear eigenvalue optimization problem can be efficiently globally solved. Moreover, we propose a convex semidefinite programming (CSDP) relaxation, which is conjectured to be the tightest among all existing upper bounds as we are unable to find a counterexample.

%



\end{document}